\documentclass[11pt]{amsart}
\usepackage{amsmath,amsthm,amssymb, amscd, amsfonts}
\usepackage[all]{xy}
\usepackage{leftidx}
\usepackage[latin1]{inputenc}
\usepackage{enumerate}
\usepackage[dvipsnames]{xcolor}
\usepackage{makecell}

\theoremstyle{plain}

\newtheorem{theorem}{Theorem}[section]
\newtheorem{corollary}[theorem]{Corollary}

\newtheorem{lemma}[theorem]{Lemma}

\theoremstyle{definition}

\theoremstyle{remark}
\newtheorem{remark}[theorem]{Remark}

\numberwithin{equation}{section}\theoremstyle{plain}

\renewcommand{\1}{\textbf{1}}

\newcommand{\C}{{\mathcal C}}
\newcommand{\D}{{\mathcal D}}
\newcommand{\F}{{\mathcal F}}

\newcommand{\Z}{{\mathcal Z}}

\newcommand{\E}{{\mathcal E}}
\newcommand{\U}{{\mathcal U}}

\newcommand{\Rep}{\operatorname{Rep}}

\newcommand\Irr{\operatorname{Irr}}
\newcommand\FPdim{\operatorname{FPdim}}

\newcommand\vect{\operatorname{Vec}}

\newcommand\id{\operatorname{id}}
\newcommand\ad{\operatorname{ad}}

\newcommand\Fun{\operatorname{Fun}}

\newtheorem*{theorem1}{Theorem \ref{mainthm}}

\begin{document}
\title[Maximally non-self-dual modular categories]{Classification of maximally non-self-dual modular categories of small dimension}
\author[F. Xu]{Fengshuo Xu$^1$}
\email{3528935465@qq.com}

\author[J. Dong]{Jingcheng Dong$^{1,2}$}
\email{jcdong@nuist.edu.cn}
\address{1. College of Mathematics and Statistics, Nanjing University of Information Science and Technology, Nanjing 210044, China}
\address{2. Center for Applied Mathematics of Jiangsu Province, Nanjing University of Information Science and Technology, Nanjing 210044, China}



\date{\today}

\begin{abstract}
 We prove that a non-pointed maximally non-self-dual (MNSD) modular category of Frobenius-Perron (FP) dimension less than $2025$ has at most two possible types, and all these types can be realized except those of FP dimension $675$, $729$ and $1125$. We also prove that all these modular categories are group-theoretical except the modular categories of dimension $675$. Our result shows that a non-group-theoretical MNSD modular category of smallest FP dimension may be the category of FP dimension $675$, and non-pointed MNSD modular category of smallest FP dimension is the category of FP dimension $243$.
\end{abstract}

\maketitle

\section{Introduction}\label{sec1}
The problem of classifying modular categories dates back the classification problem for rational conformal field theories around 1980s. Recently, modular categories with integral FP dimension have been extensively studied, for example \cite{Bruillard20162364,2012BrRoModular,2016DongNatale,HONG20101000}.

A maximally non-self-dual (MNSD) fusion category is a category that the unit object $\1$ is the only simple non-self-dual object. In this paper, we will focus our attention on MNSD modular categories and study two questions arised in \cite{2012BrRoModular}.

\medbreak
In \cite{2012BrRoModular}, Bruillard and Rowell found an example of an MNSD modular category of
rank $25$ and FP dimension $441$ that is not pointed \cite[Remark 4.6]{2012BrRoModular}. Recall
that a fusion category is called pointed if all its simple objects have FP dimension $1$.  A natural follow-up question is if there exists a non-pointed MNSD modular category of rank less than $25$ or dimension less than $441$. Czenky, Gvozdjak and Plavnik \cite{czenky2023classification} studied MNSD modular categories of rank less than $73$ and gave a positive answer to this question. Specifically, they proved that MNSD modular categories of rank less than $25$ are indeed pointed. In the present work, our first aim is to study if there exists a non-pointed MNSD modular category of FP dimension less than $441$.

Bruillard and Rowell also asked if MNSD modular categories are necessarily group-theoretical \cite[Question 4.7]{2012BrRoModular}. Jordan and Larson \cite{jordan2009classification}  constructed explicit examples of non-group-theoretical integral fusion categories. If these fusion categories have odd dimension then their Drinfeld centers are odd-dimensional non-group-theoretical modular categories. Thanks to \cite[Corollary 8.2]{Ng200734} and \cite[Theorem 2.2]{HONG20101000}, a modular category is MNSD if and only if it is odd-dimensional. Therefore, these Drinfeld centers are  examples of  non-group-theoretical  MNSD modular categories. In the present work, our second  aim is to determine the smallest dimension of which an MNSD modular category is non-group-theoretical.

\medbreak
Let $1 = d_1 < d_2< \cdots < d_s$ be positive real numbers, and let $n_1,n_2,\cdots,n_s$ be positive integers. Recall that we say a fusion category $\C$ admits type $(d_1,n_1;d_2,n_2;\cdots;d_s,n_s)$ if, for all $i = 1, \cdots, s$, $n_i$ is the number of the non-isomorphic simple objects of FP dimension $d_i$.

Much information of a fusion category is encoded in its type. For example, the FP dimension of $\C$ is $\sum_{i=1}^sn_id_i^2$ and the FP dimension of the largest pointed fusion subcategory $\C_{pt}$ is $n_1$. Conversely, if $(d_1,n_1;d_2,n_2;\cdots;d_s,n_s)$ is the type of a fusion category then it should obey strict restrictions. Our work starts from seeking such restrictions and then applies them to MNSD modular categories of FP dimension less than $2025$. Our result shows that non-pointed MNSD modular categories of FP dimension less than $2025$ are those of FP dimension $243$, $441$, $729$, $1125$, $1215$, $1323$, $1521$ and $1701$. Specifically, we prove the following theorem.

\begin{theorem1}
If $\C$ is a non-pointed MNSD modular categories of FP dimension less than $2025$ then we have the following.
\begin{itemize}
 \item If $\FPdim(\C)=243$ then $\C$ admits type $(1,9;3,26)$ or $(1,243)$. If $\C$ is not pointed then it is the modular subcategory of the Drinfeld center $\Z(\vect_{H_3}^{\omega})$, where $H_3$ denotes the Heisenberg group of order $3^3$.

 \item If $\FPdim(\C)=441$ then $\C$ admits type $(1,3;3,16;7,6)$ or $(1,441)$. If $\C$ is not pointed then it is equivalent to $\Rep(D^{\omega}(\mathbb{Z}_7\rtimes\mathbb{Z}_3))$.
     
 \item If $\FPdim(\C)=675$ then $\C$ admits type $(1,9;3,24;5,18)$ or $(1,675)$.

 \item If $\FPdim(\C)=729$ then $\C$ admits type $(1,9;3,80)$,$(1,27;3,78)$ or $(1,729)$

 \item If $\FPdim(\C)=1125$ then $\C$ admits type $(1,15;3,40;5,30)$ or $(1,1125)$.

 \item If $\FPdim(\C)=1215$ then $\C$ admits type  $(1,45;3,130)$ or $(1,1215)$. If $\C$ is not pointed then it is realized as $\D \boxtimes \E$, where $\D$ is the non-pointed modular subcategory of FP dimension $243$, $\E$ is a modular category of dimension $5$.

 \item If $\FPdim(\C)=1323$ then $\C$ admits type $(1,9;3,48;7,18)$ or $(1,1323)$. If $\C$ is not pointed then it is realized as $\Rep(D^{\omega}(\mathbb{Z}_7\rtimes\mathbb{Z}_3))\boxtimes \E$, where $\E$ is a modular category of dimension $3$.

 \item If $\FPdim(\C)=1521$ then $\C$ admits type $(1,3;3,56;13,6)$ or $(1,1521)$. If $\C$ is not pointed then it  is realized as $\Rep(D^{\omega}(\mathbb{Z}_{13}\rtimes\mathbb{Z}_3))$.

 \item If $\FPdim(\C)=1701$ then $\C$ admits type $(1,63;3,182)$ or $(1,1701)$. If $\C$ is not pointed then it is realized as $\D \boxtimes \E$, where $\D$ is the non-pointed modular subcategory of FP dimension $243$, $\E$ is a modular category of dimension $7$.
\end{itemize}
\end{theorem1}
Based on the types of MNSD modular categories of FP dimension $1125$ and $1323$, we prove that such modular categories are group-theoretical. Together with the known results in the literature, we get MNSD modular categories of FP dimension less than $2025$ are group-theoretical.

\medbreak
The organization of the paper is as follows. In Section \ref{sec2}, we recall some  basic results and notions of fusion categories. In Section \ref{sec3}, we study the restrictions that the type of a fusion category should obey. We are chiefly concerned with the type of an MNSD modular category.  In Section \ref{sec4}, we apply the results obtained in previous section to MNSD modular categories of FP dimension less than $2025$ and get our main results.

\medbreak
The reference for the theory of fusion categories are \cite{egno2015} and \cite{drinfeld2010braided}. Throughout this paper, we will work over an algebraically closed field $k$ of characteristic $0$.

\section{Preliminaries}\label{sec2}
\subsection{Group gradings on fusion categories}
Let $\C$ be a fusion category and $\Irr(\C)$ be the set of isomorphism classes of
simple objects of $\C$. We shall use $\FPdim(X)$ to denote the FP dimension of an object $X$. The FP dimension of $\C$ is defined to be the number $\FPdim (\C) =\sum_{X \in \Irr (\C)} (\FPdim X)^2$.

A fusion category $\C$ is called weakly integral if $\FPdim (\C)$ is an integer. It is called integral if $\FPdim (X)$ is an integer for every object $X$ of $\C$. By \cite[Corollary 3.11]{gelaki2008nilpotent}, an odd-dimensional fusion category is automatically integral.

A simple object $X$ is invertible if and only if its FP dimension is $1$. The set of isomorphism classes of invertible simple objects of $\C$ is denote by $G(\C)$. Then $G(\C)$ is a group with the multiplication given by tensor product. The group $G(\C)$ acts on the set $\Irr(\C)$ by left (right) tensor product. This action preserves FP dimensions. For $X\in \Irr(\C)$, $G[X]$ will denote the stabilizer of $X$ in $G(\C)$.

A fusion category is called pointed if every simple object has FP dimension $1$. If $\C$ is a pointed fusion category, then $\C$ is equivalent to the category $\vect_G^{\omega}$  of $G$-graded vector spaces with associativity constraint given by a $3$-cocycle $\omega\in H^3(G,k^{\times})$ \cite{DiPaR1991}.

\medbreak
A fusion subcategory of $\C$ is a full tensor subcategory $\D$ such that if $X$ is an object of $\C$ isomorphic to a direct summand of an object of $\D$, then $X$ lies in $\D$. If $\D$ is a fusion subcategory of $\C$ then the quotient $\FPdim(\C)/\FPdim(\D)$ is an algebraic integer \cite[Proposition 8.15]{etingof2005fusion}. We will denote $\C_{pt}$ the unique largest pointed fusion subcategory of $\C$.

\medbreak

Let $G$ be a finite group. A $G$-graded fusion category is a fusion category $\C$ admitting a decomposition $\C = \oplus_{g \in G} \C_g$ of abelian subcategories such that $\otimes:\C_g \times \C_h \to \C_{gh}$ and $(\C_g)^* = \C_{g^{-1}}$. If $\C_g\neq 0$ for all $g\in G$, we shall say this grading is faithful. In this case, we say $\C$ is a $G$-extension of the trivial component $\C_e$. If $\C$ is a $G$-extension of $\D$ then the FP dimensions of $\C_g$ are all equal for all $g\in G$, and $\FPdim (\C)=|G|\FPdim(\D)$ \cite[Proposition 8.20]{etingof2005fusion}.

The adjoint subcategory $\C_{ad}$ is the fusion subcategory of $\C$ generated by all simple objects in the decomposition $X\otimes X^*$ for all $X\in \Irr(\C)$. Recall from \cite{gelaki2008nilpotent} that a fusion category $\C$ admits a canonical faithful grading $\C = \oplus_{g\in U(\C)}\C_g$, whose trivial component $\C_e$ equivalent to the adjoint subcategory $\C_{\ad}$. Such a grading is called the universal grading of $\C$.

\medbreak
Recall that a fusion category $\C$ is said to be nilpotent if there is a sequence of
fusion categories $\C_0=\vect$, $\C_1$,$\cdots$, $\C_n=\C$ and a sequence $G_1,\cdots,G_n$ of finite groups such that $\C_i$ is obtained from $\C_{i-1}$ by a $G_i$-extension. Furthermore, $\C$ is called cyclically nilpotent if the groups $G_i$ can be chosen to be cyclic of prime order. It is well-known that fusion categories of FP dimension a power of a prime number are nilpotent by \cite[Theorem 8.28]{gelaki2008nilpotent}.

A fusion category is called weakly group-theoretical if it is Morita equivalent to a nilpotent fusion category. A fusion category is called solvable if it is Morita equivalent to a cyclically nilpotent fusion category. A fusion category is called group-theoretical if it is Morita equivalent to a pointed fusion category.

\subsection{M\"uger centralizers}
Let $\C$ be a braided fusion category. The M\"uger centralizer $\D'$ of a braided fusion subcategory $\D \subseteq \C$ is given by:
\begin{equation*}
\D' = \{ X \in \C \; | \; c_{Y, X} \circ c_{X, Y} = \id_{X,Y} \; \text{for all} \; Y \in \D \},
\end{equation*}
where $c_{X,Y}: X\otimes Y\to Y\otimes X$ is the associated braiding. The M\"uger center $\mathcal{Z}_2(\C)$ is the centralizer $\C'$ of the category itself.  The fusion category $\C$ is called non-degenerate if $\mathcal{Z}_2(\C)\simeq \vect$, and it is called symmetric if $\mathcal{Z}_2(\C)\simeq \C$, where $\vect$ is the trivial category.

A non-degenerate braided fusion category with a ribbon structure is modular. If $\C$ is modular then we have $\C_{pt} = \C_{ad}'$ \cite[Corollary 6.9]{gelaki2008nilpotent}. Thus the M\"uger center $\mathcal{Z}_2(\C_{ad})$ of $\C_{ad}$ is $\C_{ad}\cap \C_{ad}'=(\C_{ad})_{pt}$. Hence we have the lemma below.

\begin{lemma}\label{lem41}
Let $\C$ be a modular category. Then $(\C_{ad})_{pt}=\mathcal{Z}_2(\C_{ad})$ is a symmetric subcategory.
\end{lemma}

Suppose that $\C$ is modular and that $\D \subseteq \C$ is a fusion subcategory. Then $\D''=\D$, and by \cite[Theorem 3.2]{muger2003structure}, we have the identity:
\begin{equation}\label{centfactor}
\FPdim(\D) \FPdim(\D') = \FPdim(\C).
\end{equation}

Suppose that $\C$ is modular and $\C = \oplus_{g\in U(\C)}\C_g$ is the universal grading of $\C$. Then $\U(\C)\cong G(\C)$ by \cite[Theorem 6.2]{gelaki2008nilpotent}, and hence we have

\begin{equation}\label{dim}
\FPdim(\C_g) = \frac{\FPdim(\C)}{|G(\C)|}.
\end{equation}

The following theorem is known as the M\"{u}ger Decomposition Theorem  \cite[Theorem 4.2]{muger2003structure}, and will be used frequently.

\begin{theorem}\label{MugerThm}
Let $\C$ be a modular category and let $\D$ be a modular subcategory of $\C$. Then $\C$ is braided equivalent to $\D\boxtimes \D'$, where $\D'$ is the centralizer of $\D$ in $\C$. In particular, $\D'$ is also modular.
\end{theorem}

A Tannakian category is a symmetric fusion category such that it is equivalent to the representation category of a finite group as a symmetric category. The following lemma is taken from \cite[Corollary 2.50]{drinfeld2010braided}.

\begin{lemma}\label{lem25}
Let $\C$ be a symmetric fusion category. Then either $\C$ is Tannakian or $\C$ contains a Tannakian subcategory of FP dimension $\FPdim (\C)/2$. In particular, if $\FPdim (\C)$ is odd then $\C$ is Tannakian.
\end{lemma}

\subsection{Equivariantizations and de-equivariantizations}
Let $\C$ be a fusion category and let $G$ be a group acting on $\C$ by tensor autoequivalences. Recall from \cite[Section 4]{drinfeld2010braided} that we have the fusion category $\C^G$ of $G$-equivalent objects of $\C$  which is called the $G$-equivariantization of $\mathcal{\C}$.

Let $\C$ be  a fusion category and let $\mathcal{Z}(\C)$ be the Drinfeld center of $\C$. Assume that $\Rep(G)\subseteq \mathcal{Z}(\C)$ is a Tannakian subcategory which embeds into $\C$ via the forgetful functor $\mathcal{Z}(\C)\to \C$. Let $A=\Fun(G)$ be the algebra of function on $G$. It is a commutative algebra in $\mathcal{Z}(\C)$ under the embedding  above. Let $\C_G$ be the category of left $A$-modules in $\C$. It is a fusion category which is called the de-equivariantization of $\C$ by $\Rep(G)$.

There are canonical equivalences $(\C_G)^G\cong \C$  and $(\C^G)_G\cong \C$. Moreover, we have $$\FPdim(\C^G)=|G|\FPdim(\C) \mbox{\,\,and\,\,} \FPdim(\C_G)=\frac{\FPdim(\C)}{|G|}.$$

\medbreak
Let $\C^G$ be the equivariantization of $\C$ under the action $T:\underline{G}\to {\rm\underline{Aut}}_{\otimes}\C$, and let $F:\C^G\to \C$ denote the forgetful functor. For every simple object $X$ of $\C$, let $G_X = \{g\in G : T_g(X)\cong X\}$ be the inertia subgroup of $X$. Since $X$ is a simple object, there is a $2$-cocycle $\alpha_X: G_X \times G_X\to k^*$. The simple objects of $\C^G$ are parameterized by pairs $(X,\pi)$, where $X$ runs over the $G$-orbits on $\Irr(\C)$ and $\pi$ is an equivalence class of an irreducible $\alpha_X$-projective representation of $G_X$. If we use $S_{X,\pi}$ to indicate the isomorphism class of the simple object corresponding to the pair $(X,\pi)$, we get the following relation \cite[Corollary 2.13]{burciu2013fusion}:
\begin{equation}\label{eq7}
\begin{split}
\FPdim S_{X,\pi} ={\rm dim}\pi [G:G_X] \FPdim X.
\end{split}
\end{equation}

The following lemma will be used later.
\begin{lemma}\cite[Lemma 7.2]{natale2014graphs}\label{lem51}
Let $F:\C^G\to \C$ be the forgetful functor. Let $X$ be a simple object of $\C^G$. If $\FPdim(X)$ is relatively prime to the order of $G$ then $F(X)$ is a simple object of $\C$.
\end{lemma}

Let $\C$ be a modular category and $\E\cong \Rep(G)$ be a Tannakian subcategory of $\C$.  The de-equivariantization $C_G$ is a braided $G$-crossed fusion category of FP dimension $\FPdim(C)/|G|$, see \cite[Proposition 4.26]{drinfeld2010braided}. Since $\C$ is modular, the associated $G$-grading of $\C_G$ is faithful, and the trivial component $\C^0_G$ is a modular category of FP dimension $\FPdim(\C)/|G|^2$, \cite[Proposition 4.56]{drinfeld2010braided}. In particular, we have an equivalence of braided fusion categories by \cite[Corollary 3.30]{DMNO2013},
 \begin{equation}\label{equivalence}
\begin{split}
\C\boxtimes (\C_G^0)^{rev}\cong \mathcal{Z}(\C_G).
\end{split}
\end{equation}

\begin{remark}\label{rem1}
Let $\C$ be a integral modular category and $\E\cong \Rep(G)$ be a Tannakian subcategory of $\C$. Then $\C^0_G$ is also integral \cite[Corollary 4.27]{drinfeld2010braided}. Hence the exposition above shows that $\FPdim(\C)/|G|^2$ is an integer.
\end{remark}

\section{Arithmetic properties of the type of a fusion category}\label{sec3}
In the context below, we will use $\Irr_{a}(\C)$ to denote the set of isomorphism classes of simple objects of FP dimension $a$.

\begin{lemma}\label{lem2}
Let $\C$ be an integral MNSD fusion category such that it admits type $(d_1,n_1;\cdots;d_s,n_s)$. Then $d_i$ is odd and $n_i$ is even, for all $2\leq i\leq s$. In particular, $\C$ is odd-dimensional.
\end{lemma}
\begin{proof}
Since the trivial simple object $\1$ is the only self-dual simple object, the nontrivial simple objects appear in pairs such that we can write $\Irr(\C)=\{\1,X_1,\cdots, X_t, X_1^*,\cdots,X_t^*\}$. It follows that $n_i$ is even for all $2\leq i\leq s$.

Assume on the contrary that there exists a simple object $X$ such that it has even dimension. We may write the decomposition of $X\otimes X^*$ as follows:
$$X\otimes X^*=\1\oplus \sum_{i=1}^{t}m(X_i,X\otimes X^*)X_i\oplus \sum_{i=1}^{t}m(X_i^*,X\otimes X^*)X_i^*,$$
where $m(X_i,X\otimes X^*)=m(X_i^*,X\otimes X^*)$. Counting FP dimension on both sides, we have
$$\FPdim (X)^2=1+2 \sum_{i=1}^{t}m(X_i,X\otimes X^*)\FPdim (X_i).$$
This is a contradiction since the left side is even and the right side is odd. This proves every $d_i$ is odd for all $2\leq i\leq s$.
\end{proof}

\begin{remark}
A similar result in the Hopf algebra setting was obtained in \cite{kashina2002self}.
\end{remark}

The lemma below was initially obtained in the Hopf setting, see \cite[Lemma 2.2 and Lemma 2.5]{DongDai2012}. In fact, the proof in \cite{DongDai2012} only uses the properties of the Grothendieck ring of a semisimple Hopf algebra. Therefore, the proof of \cite{DongDai2012} also works in the integral fusion category setting.
\begin{lemma}\label{lem-dim3}
Let $\C$ be an integral MNSD fusion category  such that it admits type $(d_1,n_1;\cdots;d_s,n_s)$. Assume that $d_2=3$. Then $3$ divides the order of the stabilizer $G[X]$ for every simple object $X$ in $\Irr_3(\C)$. Furthermore, if $\C$ does not have simple objects of dimension $9$ then simple objects of dimension $1$ and $3$ generate a fusion subcategory of $\C$, and hence $n_1+9n_2$ divides $\FPdim(\C)$.
\end{lemma}

\begin{lemma}\label{lem1}
Let $\C$ be a integral fusion category such that it admits type $(d_1,n_1;\cdots;d_s,n_s)$. Assume that $n_1=p^tm$, where $p$ is a prime not dividing $m$. If there exists $n_i$ $(2\leq i\leq s)$ not divisible by $p$ then there exists a simple object $X$ of dimension $d_i$ such that $p^t$ divides the order of the stabilizer $G[X]$ of $X$. In particular, if $\C$ is an MNSD modular tensor category then $p^{2t}$ divides $\FPdim(\C)$.
\end{lemma}
\begin{proof}
Let $G$ be a subgroup of $G(\C)$ of order $p^t$ and consider the action of $G$ on $\Irr_{d_i}(\C)$. The set $\Irr_{d_i}(\C)$ thus  is a union of orbits which have length $1$ or $p^i$, $1\leq i\leq t$. Since $p$ does not divide $n_i=|\Irr_{d_i}(\C)|$, there exists at least one orbit with length $1$. That is, there exists a simple object $X\in \Irr_{d_i}(\C)$ such that $g\otimes X=X$ for all $g\in G$. Hence $G$ is a subgroup of the stabilizer $G[X]$ of $X$. Thus $p^t$ divides $|G[X]|$.

From the arguments above, we know that the pointed fusion subcategory $\D$ generated by $G$ is contained in the adjoint subcategory $\C_{ad}$. If $\C$ is an MNSD modular tensor category then $\D\cong \Rep(G)$ is a Tannakian subcategory by Lemma \ref{lem41} and Lemma \ref{lem25}. Let $\D_G$ be the de-equivariantization of $\C$ by $\Rep(G)$. Then $\D_G=\oplus_{g\in G}(\D_G)_g$ has a faithful $G$-grading by \cite[Proposition 4.56]{drinfeld2010braided}. Moreover, $\D_G$ and the trivial componen $\D_G^0$ are integral fusion categories by \cite[Corollary 4.27]{drinfeld2010braided}. Hence $\FPdim (\D_G^0)=\FPdim (\C)/p^{2t}$ is an integer. This proves the result.
\end{proof}

\begin{lemma}\label{lem001}
Let $\C$ be a modular category such that it admits type\\ $(d_1,n_1;\cdots;d_s,n_s)$. Then $n_1d_i^2\leq\FPdim(\C)$ for all $i=2,\cdots,s$.
\end{lemma}
\begin{proof}
 Considering the universal grading of $\C$, we get $\U(\C)\cong G(\C)$ by \cite[Theorem 6.2]{gelaki2008nilpotent}. Hence $|\U(\C)|=n_1$. Since every component of the universal grading has the same dimension, every component has dimension at least $d_i^2$ for all $i=2,\cdots, s$. Hence the modular category $\C$ has dimension at least $n_1d_i^2$.
\end{proof}

\begin{lemma}\label{lem3}
Let $\C$ be a modular category such that it admits type $(1,n;d,m)$. Then $n$ divides $m$ or $d^2n$ divides $\FPdim(\C)$.
\end{lemma}
\begin{proof}
Let $\C=\oplus_{g\in \U(\C)}\C_g$ be the universal grading of $\C$. Then $|\U(\C)|=n$. It follows that the component $\C_g$ has dimension $1+d^2m/n$ for every $g\in \U(\C)$.

Let $a_g$ and $b_g$ be the number of non-isomorphic invertible objects and $d$-dimensional simple objects in $\C_g$, respectively. Then
$$a_g+d^2b_g=\frac{n+d^2m}{n}.$$
It is clear that there exists at least one component $\C_g$ such that $a_g=0$ or $1$, otherwise the total number of non-isomorphic invertible simple objects is at least $2|\U(\C)|=2n$, a contradiction.

If $a_g=0$ then $d^2b_g=1+d^2m/n$ which implies $(d^2n)b_g=n+d^2m=\FPdim(\C)$, hence $d^2n$ divides $\FPdim(\C)$.

If $a_g=1$ then $d^2(b_g-m/n)=0$ which implies $b_g=m/n$ is an integer, hence $n$ divides $m$.
\end{proof}

Let $\C$ be a modular category such that it admits type $(1,n;d,m)$. Then Lemma \ref{lem3} shows that $d^2n$ divides $\FPdim(\C)$. The lemma below shows that, in some situation, the extreme case $d^2n=\FPdim(\C)$ can not hold.

\begin{lemma}\label{lem3-1}
Let $\C$ be a modular category such that it admits type $(1,n;d,m)$. Assume that $\C$ is integral. Then $d^2n<\FPdim(\C)$.
\end{lemma}
\begin{proof}
Assume on the contrary that $d^2n=\FPdim(\C)$. Then every component $\C_g$ in the universal grading $\C=\oplus_{g\in \U(\C)}\C_g$ has dimension $d^2$. Hence every component $\C_g$ admits type $(1,d^2)$ or $(d,1)$. In particular, $\C_{ad}$ admits type $(1,d^2)$. By Lemma \ref{lem41}, $\C_{ad}=(\C_{ad})_{pt}$ is a symmetric subcategory of $\C$. Since $\FPdim(\C_{ad})=d^2\geq 4$, $\C_{ad}$ contains a nontrivial Tannakian subcategory by Lemma \ref{lem25}. That is, there exists $g\in G(\C_{ad})$ such that $\theta_g=1$, where $\theta$ is the ribbon structure of $\C$, see \cite{CP2022} and reference therein. On the other hand, the type $(d,1)$ implies that $g\otimes X=X$ for every simple object $X$ in $\Irr_d(\C)$. It follows from \cite[Lemma 2.3]{CP2022} that the rows of the $S$-matrix corresponding to (the isomorphism classes of ) $g$ and $1$ are equal. This is impossible since $S$ is invertible.
\end{proof}

In the lemma below, we collet two known results which will be used later.

\begin{lemma}\label{lem002}
Let $\C$ be an integral modular category such that it admits type $(d_1,n_1;\cdots;d_s,n_s)$. Then

(1)\, $d_i^2$ divides $\FPdim(\C)$ for all $i=2,\cdots,s$.

(2)\, $n_1$ divides $n_id_i^2$ for all $i=2,\cdots,s$.
\end{lemma}
\begin{proof}
Parts (1) and (2) follows from \cite[Theorem 2.11]{etingof2011weakly} and \cite[Lemma 2.2]{dong2012frobenius}, respectively.
\end{proof}

\begin{lemma}\label{prop1}
Let $\C$ be a modular category of FP dimension $p^3q^sm$ or $p^tm$, where $p,q$ are distinct odd prime numbers, $m$ is a square-free integer with $(m,pq)=1$, and $s\leq 4$, $t\leq 6$. Assume that $\C$ admits type $(1,p;p,n_2;q,n_3;\cdots)$, where $(p,n_2)=1$. Then the adjoint subcategory $\C_{ad}$ can not admit type $(1,p;p,n_2';q,n_3';\cdots)$.
\end{lemma}
\begin{proof}
Considering the action of $G(\C)$ on $\Irr_{p}(\C)$, we know $\Irr_{p}(\C)$ is a union of orbits which have length $1$ or $p$. Since $(p,n_2)=1$, there exists at least one orbit with length $1$. Hence there exists a simple object $X$ of dimension $p$ such that $G[X]=G(\C)$. This means that the pointed fusion subcategory $\C_{pt}$ is contained in $\C_{ad}$. By Lemma \ref{lem41}, $(\C_{ad})_{pt}=\C_{pt}$ is a symmetric subcategory. Since $\FPdim(\C_{pt})=p$ is odd, $\C_{pt}\cong\Rep(\mathbb{Z}_p)$ is a Tannakian subcategory.

Assume  on the contrary that $\C_{ad}$ admits the type $(1,p;p,n_2';q,n_3';\cdots)$. Let $(\C_{ad})_{\mathbb{Z}_p}$  be the de-equivariantization of $\C_{ad}$ by $\mathbb{Z}_p$. Then $\FPdim((\C_{ad})_{\mathbb{Z}_p})=pq^sm$ or $p^{t-2}m$. In addition, $(\C_{ad})_{\mathbb{Z}_p}$ is a modular category by \cite[Remark 2.3]{etingof2011weakly}. It follows from \cite[Corollary 4.13]{2016DongNatale} that $(\C_{ad})_{\mathbb{Z}_p}$ is pointed. But Lemma \ref{lem51} shows that the image $F(X)$ of any simple object $X$ of dimension $q$ under the functor $F:\C_{ad}\to (\C_{ad})_{\mathbb{Z}_p}$ is again a simple object in $(\C_{ad})_{\mathbb{Z}_p}$. This is a contradiction.
\end{proof}

\section{MNSD modular categories of dimension $<2025$}\label{sec4}
In this section, we shall determine all possible types of MNSD modular categories of dimension less than $2025$. As we mentioned in the Introduction, these modular categories are odd-dimensional. By \cite[Theorem 8.2]{natale2013weakly}, all these modular categories are solvable. Hence if they have types $(d_1,n_1;\cdots;d_s,n_s)$ then $n_1>1$. In addition, \cite[Corollary 4.3]{2016DongNatale} shows that modular categories of dimension $dq^n$ is pointed, where $q$ is an odd prime number, $d$ is a square-free integer not divisible by $q$ and $n\leq 4$. Hence it suffices to consider the cases $\FPdim(\C)=225=3^2\times 5^2$, $243=3^5$, $441=3^2\times 7^2$, $675=3^3\times 5^2$, $729=3^6$, $1089=3^2\times 11^2$,$1125=3^2\times 5^3$, $1215=3^5\times 5$, $1225=5^2\times 7^2$,$1323=3^3\times 7^2$, $1521=3^2\times 13^2$, $1575=3^2\times 5^2\times 7$ and $1701=3^5\times 7$.

\begin{theorem}\label{mainthm}
Let $\C$ be a modular categories of dimension $225$, $243$, $441$, $675$, $729$, $1089$, $1125$, $1215$, $1225$, $1323$, $1521$, $1575$ or $1701$. Then we have:
\begin{itemize}
 \item If $\FPdim(\C)=225$, $1089$, $1225$ or $1575$ then $\C$ is pointed.

 \item If $\FPdim(\C)=243$ then $\C$ admits type $(1,9;3,26)$ or $(1,243)$. If $\C$ is not pointed then it is the modular subcategory of the Drinfeld center $\Z(\vect_{H_3}^{\omega})$, where $H_3$ denotes the Heisenberg group of order $3^3$.

 \item If $\FPdim(\C)=441$ then $\C$ admits type $(1,3;3,16;7,6)$ or $(1,441)$. If $\C$ is not pointed then it is equivalent to $\Rep(D^{\omega}(\mathbb{Z}_7\rtimes\mathbb{Z}_3))$.
     
 \item If $\FPdim(\C)=675$ then $\C$ admits type $(1,9;3,24;5,18)$ or $(1,675)$.

 \item If $\FPdim(\C)=729$ then $\C$ admits type $(1,9;3,80)$,$(1,27;3,78)$ or $(1,729)$.

 \item If $\FPdim(\C)=1125$ then $\C$ admits type $(1,15;3,40;5,30)$ or $(1,1125)$.

 \item If $\FPdim(\C)=1215$ then $\C$ admits type  $(1,45;3,130)$ or $(1,1215)$. If $\C$ is not pointed then it is realized as $\D \boxtimes \E$, where $\D$ is the non-pointed modular subcategory of FP dimension $243$, $\E$ is a modular category of dimension $5$.

 \item If $\FPdim(\C)=1323$ then $\C$ admits type $(1,9;3,48;7,18)$ or $(1,1323)$. If $\C$ is not pointed then it is realized as $\Rep(D^{\omega}(\mathbb{Z}_7\rtimes\mathbb{Z}_3))\boxtimes \E$, where $\E$ is a modular category of dimension $3$.

 \item If $\FPdim(\C)=1521$ then $\C$ admits type $(1,3;3,56;13,6)$ or $(1,1521)$. If $\C$ is not pointed then it  is realized as $\Rep(D^{\omega}(\mathbb{Z}_{13}\rtimes\mathbb{Z}_3))$.

 \item If $\FPdim(\C)=1701$ then $\C$ admits type $(1,63;3,182)$ or $(1,1701)$. If $\C$ is not pointed then it is realized as $\D \boxtimes \E$, where $\D$ is the non-pointed modular subcategory of FP dimension $243$, $\E$ is a modular category of dimension $7$.
\end{itemize}
%
%
%
%
%
%
%
%
%
%
%
%
\end{theorem}

We write a Matlab program by which we first find all possible positive integers $1=d_1<d_2<\cdots<d_s$
and $n_1,n_2,\cdots,n_s$ such that $\FPdim(\C)=\sum_{i=1}^{s}n_id_i^2$, and then exclude those which are not possible types by using Lemmas \ref{lem2}-\ref{lem002}. We hence get the following possible types.

$\FPdim(\C)=225: (1,3;3,8;5,6),(1,225)$.

$\FPdim(\C)=243: (1,9;3,26),(1,243)$.

$\FPdim(\C)=441: (1,3;3,16;7,6),(1,441)$.

$\FPdim(\C)=675:(1,3;3,8;5,6;15,2),	(1,3;3,8;5,24),	(1,9;3,24;5,18)$,\\ $(1,675)$.

$\FPdim(\C)=729: (1,9;3,8;9,8),	(1,9;3,26;9,6),	(1,9;3,80),(1,27;3,78)$,\\ $(1,729)$.

$\FPdim(\C)=1089: (1,3;3,40;11,6), (1,1089)$.

$\FPdim(\C)=1125: (1,3;3,8;5,6;15,4), (1,3;3,8;5,42)$, $(1,3;3,8;5,24;15,2)$, $(1,9;3,24;5,36)$,	$(1,15;3,40;5,30)$, $(1,1125)$.

$\FPdim(\C)=1215: (1,9;3,8;9,14)$,	$(1,9;3,26;9,12)$,	$(1,9;3,44;9,10)$, $ (1,9;3,134),	(1,27;3,132)$,	$(1,45;3,130)$, $(1,1215)$.

$\FPdim(\C)=1225: (1,1225)$.

$\FPdim(\C)=1323: (1,3;3,16;7,6;21,2), (1,3;3,16;7,24), (1,9;3,48;7,18)$, $(1,21;3,14;7,24)$, $(1,1323)$.

$\FPdim(\C)=1521: (1,3;3,56;13,6),	(1,1521)$.

$\FPdim(\C)=1575: $ $(1,3;3,8;5,6;15,6)$,	$(1,3;3,8;5,24;15,4)$,\\	$(1,3;3,8;5,42;15,2)$,	$(1,3;3,8;5,60)$,	$(1,3;3,58;5,6;15,4)$,	$(1,3;3,58;5,42)$,	$(1,9;3,24;5,54)$,$(1,3;3,58;5,24;15,2)$,	$(1,21;3,56;5,42)$,	$(1,1575)$.

$\FPdim(\C)=1701: $ $(1,9;3,8;9,20)$,	$(1,9;3,26;9,18)$,	$(1,9;3,62;9,14)$,	$(1,9;3,188)$,	$(1,27;3,186)$,	$(1,63;3,182)$,	$(1,1701)$,

\medbreak
Our proof involves Lemmas \ref{lem42}-\ref{lem49} and each lemma will deal with one or several dimensions. We will prove our statement discarding the possible types until we are left with the types stated above. In our arguments below, we will frequently use the fact $\FPdim(\C_{ad})=\FPdim(\C)/n_1$, see equation (\ref{dim}).

\begin{lemma}\label{lem42}
Let $\C$ be a modular category of FP dimension $225$, $1089$, $1225$ or $1575$. Then $\C$ is pointed.
\end{lemma}

\begin{proof}
First, a modular category of FP dimension $1225$ is pointed by the data listed above.

\medbreak
\textbf{Case $\FPdim(\C)=225$ or $1089$.} Let $p<q$ be prime number. Recall from \cite[Theorem 4.2]{bruillard2013classification}, if a $p^2q^2$-dimensional modular category $\C$ has $|G(\C)|=p$ then $p|q-1$. Hence types $(1,3;3,8;5,6)$ ($\FPdim(\C)=225$) and $(1,3;3,40;11,6)$ ($\FPdim(\C)=1089$) can be discarded. Hence they are pointed.

\medbreak
\textbf{Case $\FPdim(\C)=1575$.} Since $\FPdim(\C) = 3^2\times 5^2\times 7$, $\C$ is weakly group-theoretical by \cite[Theorem 7.4]{natale2013weakly}. It follows from  \cite[Theorem 4.2]{OSTRIK2023108961} that we have a decomposition $\C = \D\boxtimes \E$, where $\D$ and $\E$ are modular categories of FP dimension $3^2\times 5^2$ and $7$, respectively. By Lemma \ref{lem42}, $\D$ is pointed and so is $\C$.
\end{proof}

\begin{lemma}\label{lem43}
Let $\C$ be a modular category of FP dimension $243$. Then $\C$ admits type $(1,9;3,26)$ or $(1,243)$. If $\C$ is not pointed then it is the modular subcategory of the Drinfeld center $\Z(\vect_{H_3}^{\omega})$, where $H_3$ denotes the Heisenberg group of order $3^3$.
\end{lemma}
\begin{proof}
By the data listed above,  $\C$ admits type  $(1,9;3,26)$ or $(1,243)$. If $\C$ admits type $(1,9;3,26)$ then $\C$ has rank $35$. By \cite[Theorem 1.1]{czenky2023classification}, $\C$ is the modular subcategory of the Drinfeld center $\Z(\vect_{H_3}^{\omega})$ as desired.
\end{proof}

\begin{corollary}\label{cor47}
Any modular category of FP dimension $3^5d$ with $d$ square-free is either pointed, or equivalent to $\D\boxtimes \E$,  where $\D$ is the modular subcategory of the Drinfeld center $\Z(\vect_{H_3}^{\omega})$ as in Lemma \ref{lem43} and $\E$ is a pointed modular category of dimension $d$.
In particular, modular categories of FP dimension $1215$ or $1701$ have such structure.
\end{corollary}
\begin{proof}
Let $\D$ be a modular category of FP dimension $3^5d$. Then $\D$ is nilpotent by \cite[Theorem 4.7]{2016DongNatale}, and hence we have a decomposition $\D\cong \E \boxtimes \F$, where $\FPdim(\E)=3^5$ and $\FPdim(\F)=d$, see \cite[Theorem 6.10]{drinfeld2007g}. Moreover, Both $\E$ and $\F$ are modular by \cite[Remark 2.1]{2016DongNatale}, and the result then follows from Lemma \ref{lem43}. Since $1215=3^5\times 5$ and $1701=3^5\times 7$, the previous result can apply.
\end{proof}

\begin{lemma}\label{lem44}
Let $\C$ be a modular category of FP dimension $441$. Then $\C$ admits type $(1,3;3,16;7,6)$ or $(1,441)$. If $\C$ is not pointed then it is equivalent to $\Rep(D^{\omega}(\mathbb{Z}_7\rtimes\mathbb{Z}_3))$.
\end{lemma}
\begin{proof}
By the data listed above,  $\C$ admits type  $(1,3;3,16;7,6)$ or $(1,441)$. If $\C$ admits type $(1,3;3,16;7,6)$ then it has rank $25$. By \cite[Theorem 1.1]{czenky2023classification},  it is equivalent to $\Rep(D^{\omega}(\mathbb{Z}_7\rtimes\mathbb{Z}_3))$.
\end{proof}

\begin{lemma}\label{lem44-1}
Let $\C$ be a modular category of FP dimension $675$. Then $\C$ admits type $(1, 9; 3, 24; 5, 18)$ or $(1,675)$.
\end{lemma}
\begin{proof}
\medbreak
If $\C$ admits type $(1, 3; 3, 8; 5, 6; 15, 2)$ or $(1, 3; 3, 8;$ $5, 24)$ then $225=\FPdim(\C_{ad})$. Counting dimension, we get that $\C_{ad}$ admits type $(1, 3; 3, 8; 5, 6)$. By Lemma \ref{prop1}, it is impossible. Hence only types $(1, 9; 3, 24; 5, 18)$ and $(1,675)$ are left.
\end{proof}

\begin{lemma}\label{lem45}
Let $\C$ be a modular category of FP dimension $729$. Then $\C$ admits type $(1,9;3,80)$, $(1,27;3,78)$ or $(1,729)$.
\end{lemma}
\begin{proof}
If $\C$ admits type $(1,9;3,8;9,8)$ or  $(1,9;3,26;9,6)$ then  $\FPdim(\C_{ad})=81$. Counting dimension, we get that $\C_{ad}$ admits type $(1,9;3,8)$. Hence $(\C_{ad})_{pt}=\C_{pt}$ is Tannakian by Lemma \ref{lem41} and Lemma \ref{lem25}. Let $\D=\Rep(\mathbb{Z}_3)$ be a subcategory of $\C_{pt}$ of dimension $3$. Then $\FPdim(\D')=243$ by identity (\ref{centfactor}). Counting dimension, we get that $\D'$ admits type $(1,9;3,8;9,2)$. The M\"{u}ger center $\mathcal{Z}_2(\D')$ of $\D'$ is $\D''\cap \D'=\D\cap \D'=\D$. Let $\D'_{\mathbb{Z}_3}$ be the de-equivariantization of $\D'$ by $\mathbb{Z}_3$. Then $\D'_{\mathbb{Z}_3}$ is a modular category by \cite[Remark 2.3]{etingof2011weakly}. Since $\FPdim(\D'_{\mathbb{Z}_3})=81$, it is pointed by \cite[Corollary 4.3]{2016DongNatale}. It follows that the simple objects of $\D'$ has dimension at most $3$ by equation (\ref{eq7}). This contradicts the type of $\D'$. Hence types $(1,9;3,8;9,8)$ and   $(1,9;3,26;9,6)$ can be discarded. Hence only types $(1,9;3,80)$, $(1,27;3,78)$ and $(1,729)$ are left.
\end{proof} 

\begin{lemma}\label{lem46}
Let $\C$ be a modular category of FP dimension $1125$. Then $\C$ admits type $(1,15;3,40;5,30)$ or $(1,1125)$.
\end{lemma}
\begin{proof}
If $\C$ admits type $(1; 3; 3; 8; 5; 24; 15; 2)$, $(1; 3; 3; 8; 5; 6; 15; 4)$ or $(1; 3; 3; 8;$ $5; 42)$
then $\FPdim(\C_{ad}) = 375$. Counting dimension, we get that $\C_{ad}$ admits type $(1; 3; 3; 8; 5; 12)$. By Lemma \ref{prop1}, it is impossible.

If $\C$ admits type $(1; 9; 3; 24; 5; 36)$ then $\FPdim(\C_{ad}) = 125$ which shows that $\C_{ad}$ is nilpotent, and so is $\C$. By \cite[Theorem 6.10]{drinfeld2007g}, we have a decomposition $\C = \D\boxtimes \E$, where $\FPdim(\D) = 3^2$ and $\FPdim(\E) = 5^3$. Clearly,$\D$ is pointed since it can not contain a simple object of dimension 3.

Hence $\D = \C_{pt}$ is the largest pointed fusion subcategory of $\C$. It follows that the largest pointed
fusion subcategory $E_{pt}$ of $\E$ is trivial. This contradicts Lemma \ref{lem-dim3} since $\E$ contains simple objects of dimension $3$. Hence only types $(1,15;3,40;5,30)$ or $(1,1125)$ are left.
\end{proof}

\begin{lemma}\label{lem48}
Let $\C$ be a modular category of FP dimension $1323$. Then $\C$ admits type $(1,9;3,48;7,18)$ or $(1,1323)$. If $\C$ is not pointed then it can be realized as $\Rep(D^{\omega}(\mathbb{Z}_7\rtimes\mathbb{Z}_3))\boxtimes \E$, where $\E$ is a modular category of dimension $3$.
\end{lemma}
\begin{proof}
If $\C$ admits type $(1,3;3,16;7,6;21,2)$ or $(1,3;3,16;7,24)$ then $\C_{ad}$ has dimension $441$. Counting dimension, we get that $\C_{ad}$ admits type $(1,3;3,16;7,24)$. By Proposition \ref{prop1}, it is impossible.

If  $\C$ admits type $(1,21;3,14;7,24)$ then Lemma \ref{lem1} can be applied as follows: $n_1=3\times7$, $n_2=14$ is not divisible by $3$, hence $(\C_{ad})_{pt}$ has a fusion subcategory of dimension $3$; $n_2=24$ is not divisible by $7$, hence $(\C_{ad})_{pt}$ has a fusion subcategory of dimension $7$. It follows that $(\C_{ad})_{pt}=\C_{pt}$. Considering the equation below obtained from the possible type of $\C_{ad}$:
$$21+9x+49y=\FPdim(\C_{ad}),$$
It has no solutions.  Hence such adjoint subcategory can not exist.

If  $\C$ admits type $(1,27;3,144)$ then $\FPdim(\C_{ad})=49$. Hence $\C_{ad}$ is nilpotent and so is $\C$. By \cite[Corollary 5.3]{gelaki2008nilpotent}, $\FPdim(X)=1$ or $7$ for all simple object $X$ in $\C$, a contradiction.

Hence, only types $(1,9;3,48;7,18)$ and $(1,1323)$ are left. Obviously, type $(1,9;3,48;7,18)$ can be realized as $\Rep(D^{\omega}(\mathbb{Z}_7\rtimes\mathbb{Z}_3))\boxtimes \E$, where $\E$ is a modular category of dimension $3$.
\end{proof}

\begin{lemma}\label{lem49}
Let $\C$ be a modular category of FP dimension $1521$. Then $\C$ admits type $(1,3;3,56;13,6)$ or $(1,1521)$. If $\C$ is not pointed then it can be realized as $\Rep(D^{\omega}(\mathbb{Z}_{13}\rtimes\mathbb{Z}_3))$.
\end{lemma}

\begin{proof}
By the data listed above,  $\C$ admits type $(1,3;3,56;13,6)$ or $(1,1521)$. Obviously, type $(1,3;3,56;13,6)$ can be realized as $\Rep(D^{\omega}(\mathbb{Z}_{13}\rtimes\mathbb{Z}_3))$.
\end{proof}

\begin{corollary}
MNSD modular categories of FP dimension less than $2025$ are group-theoretical except the modular categories of dimension $675$.
\end{corollary}

\begin{proof}
 By \cite[Theorem 4.2]{bruillard2013classification} and \cite[Corollary 4.8]{2016DongNatale}, they are group-theoretical if they have dimension $p^2q^2$ or $p^nd$, where $p,q$ are distinct prime numbers and $d$ is a square-free integer not divisible by $p$. It suffices to consider the cases $\FPdim(\C)=1125$ and $1323$. Moreover, pointed fusion categories are automatically group-theoretical. Therefore, it is enough to consider non-pointed cases.
\medbreak
\textbf{Case $\FPdim(\C)=1125$.} Since $\C$ admits type $(1,15;3,40;5,30)$, the group $G(\C)$ has order $15$. Hence the largest pointed fusion subcategory $\C_{pt}$ has a fusion subcategory $\D$ of dimension $5$. The fusion category $\D$ is either modular or Tannakian. If $\D$ is modular then $\C\cong \D\boxtimes \E$ by the M\"{u}ger decomposition Theorem, where $\E$ has dimension $225$ and type $(1,3;3,8;5,6)$. By Lemma \ref{lem42}, $\E$ can not exist. Hence $\D$ is Tannakian.

By Lemma \ref{lem-dim3}, $(\C_{ad})_{pt}$ contains a fusion subcategory of dimension $3$. It is Tannakian by  Lemma \ref{lem41} and Lemma \ref{lem25}. Hence $\C_{pt}$ is Tannakian. Assume there exists $G$ such that $\C_{pt}\cong\Rep(G)$ and $\C_G$ be the de-equivariantization of $\C$ by $G$. Then $\C_G$ has a faithful $G$-grading $\C_G=\oplus_{g\in G}(\C_G)_g$ such that $\C_G^0$ is modular and $\FPdim(\C_G^0)=5$. Hence $\C_G$ is nilpotent since $\C_G^0$ is nilpotent. By \cite[Corollary 5.3]{gelaki2008nilpotent}, $\FPdim(X)=1$ for all simple object $X$ in $\C_G$. Thus $\C$ is group-theoretical by \cite[Theorem 7.2]{naidu2009fusion}.

\medbreak
\textbf{Case $\FPdim(\C)=1323$.} In this part, we will follow the proof of \cite[Theorem 3.19]{czenky2023classification} to get a Tannakian subcategory of dimension $21$. Since $\C$ admits type $(1,9;3,48;7,18)$,  Lemma \ref{lem-dim3} shows that $(\C_{ad})_{pt}$ contains a fusion subcategory of dimension $3$. In addition, $(\C_{ad})_{pt}\neq \C_{pt}$ otherwise $91$ divides $\FPdim(\C)$ by Remark \ref{rem1}. Hence we may assume $(\C_{ad})_{pt}\cong \Rep(\mathbb{Z}_3)$ as a Tannakian subcategory. Let $(\C_{ad})_{\mathbb{Z}_3}$ be the de-equivariantization of $\C_{ad}$ by $\mathbb{Z}_3$. Then $(\C_{ad})_{\mathbb{Z}_3}$ is modular \cite[Remark 2.3]{etingof2011weakly} and $\FPdim((\C_{ad})_{\mathbb{Z}_3})=49$ \cite[Proposition 4.26]{drinfeld2010braided}. Hence $(\C_{ad})_{\mathbb{Z}_3}$ is pointed by \cite[Corollary 4.13]{2016DongNatale}. Hence $(\C_{ad})_{\mathbb{Z}_3}\cong \vect_G^{\omega}$, where $G$ is an abelian group of order $49$.

Since $G\cong \mathbb{Z}_{49}$ or $\mathbb{Z}_{7}\times\mathbb{Z}_{7}$, it has $1$ or $8$ subgroups of order $7$, respectively. Since $3\nmid 8$, then there exists at least one subgroup of order $7$ that is invariant under the action of $G$. Let $\D$ be the fusion subcategory of $(\C_{ad})_{\mathbb{Z}_3}$ associated with said subgroup which has FP dimension $7$. Then the equivariantization $\D^{\mathbb{Z}_3}$ is a fusion subcategory of $\C_{ad}$ and has FP dimension $21$. It is not pointed since $(\C_{ad})_{pt}$ only has dimension $3$. Hence $\D^{\mathbb{Z}_3}$ admits type $(1,3;3,2)$. By the classification of premodular categories of rank $5$ \cite{BrMar2018}, $\D^{\mathbb{Z}_3}\cong\Rep(\mathbb{Z}_7\rtimes\mathbb{Z}_3)$ is a Tannakian subcategory.

Let $\C_G$ be the de-equivariantization of $\C$ by $G$, where $G=\mathbb{Z}_7\rtimes\mathbb{Z}_3$. Then $\C_G$ has a faithful $G$-grading $\C_G=\oplus_{g\in G}(\C_G)_g$ such that $\C_G^0$ is modular and $\FPdim(\C_G^0)=3$. Since $\C_G^0$ is nilpotent, $\C_G$ is nilpotent. By \cite[Corollary 5.3.]{gelaki2008nilpotent}, $\FPdim(X)=1$ for all simple object $X$ in $\C_G$. That is, $\C_G$ is pointed. It follows that $\C$ is group-theoretical by \cite[Theorem 7.2]{naidu2009fusion}.
\end{proof}

\section*{Acknowledgements}
The research of the author is supported by the Natural Science Foundation of Jiangsu Providence (Grant No. BK20201390).



\end{document}